\font\script=eusm10.
\font\sets=msbm10.
\font\stampatello=cmcsc10.
\def\0{{\bf 0}}
\def\1{{\bf 1}}
\def\defineq{\buildrel{def}\over{=}}
\def\definiz{\buildrel{def}\over{\Longleftrightarrow}}

\def\C{\hbox{\sets C}}
\def\N{\hbox{\sets N}}
\def\P{\hbox{\sets P}}

\def\cloud0{\hbox{$<\0>$}}
\def\corsivoF{\hbox{\script F}}
\def\corsivoFzero{\hbox{\script F}_0}
\def\corsivoN{\hbox{\script N}}
\def\corsivoP{\hbox{\script P}}
\def\corsivoQ{\hbox{\script Q}}
\def\corsivoR{\hbox{\script R}}
\def\corsivoS{\hbox{\script S}}

\def\Wintner{{\rm Win}}
\def\square{\hbox{\vrule\vbox{\hrule\phantom{s}\hrule}\vrule}}
\def\supporto{{\rm supp}}
\par
\centerline{\bf Finite and infinite Euler products of Ramanujan expansions}
\bigskip
\bigskip
\centerline{Giovanni Coppola}\footnote{ }{MSC$2010$:$11{\rm A}25$,$11{\rm K}65$,$11{\rm N}37$ - Keywords: Ramanujan expansion, Euler product, finite Euler product} 
\bigskip
\bigskip
\rightline{\it to Euler aficionados} 
\bigskip
\bigskip
\par
\noindent
{\bf Abstract}. All the $F:\N \rightarrow \C$ having Ramanujan expansion $F(a)=\sum_{q=1}^{\infty}G(q)c_q(a)$ (here $c_q(a)$ is the Ramanujan sum) pointwise converging in $a\in \N$, with $G:\N \rightarrow \C$ a multiplicative function, may be factored into two Ramanujan expansions, one of which is  a {\it finite Euler product} : see our Main Theorem. This is a general result, with unexpected and useful consequences, esp., for the Ramanujan expansion of null-function, say $\0$. The Main Th.m doesn't require other analytic assumptions, as pointwise convergence suffices; this depends on a general property of Euler $p-$factors (the factors in Euler products) for the general term $G(q)c_q(a)$; namely, once fixed $a\in \N$ (and prime $p$), the $p-$Euler factor of $G(q)c_q(a)$ (involving all $p-$powers) has a finite number of non-vanishing terms (depending on $a$): see our Main Lemma. In case we also add some other hypotheses, like the absolute convergence, we get more classical Euler products: the infinite ones. For the Ramanujan expansion of $\0$ this strong hypothesis makes the class of $\0$ Ramanujan coefficients much smaller; also excluding Ramanujan's \thinspace $G(q)=1/q$ and Hardy's \thinspace $G(q)=1/\varphi(q)$ ($\varphi$ is Euler's totient function). 
\bigskip
\bigskip
\bigskip
\par
\noindent {\bf 1. Introduction, statement and proof of main results}

\bigskip

\par
\noindent
In [C3] we studied the Ramanujan expansions of very general arithmetic functions $F:\N \rightarrow \C$, only assuming the Ramanujan Conjecture for their growth, namely $F(n)\ll_{\varepsilon}n^{\varepsilon}$. Here $\ll$ is classical Vinogradov notation, meaning that LHS (Left Hand Side) is bounded in modulus by a constant $C>0$ times the RHS (Right Hand Side), as $n\to \infty$; the subscript, here $\varepsilon>0$ an arbitrarily small constant, expresses an eventual dependence of $C$ from variables therein. This, of course, is not a kind of strong assumption, as we may (at least for reasonable $F$) think to re-normalize our arithmetic function. However, the real constraint we put there [C3] on $F$ is having a, say, Eratosthenes Transform [W] $F'\defineq \mu \ast F$ (see [T] for definitions and properties of {\it M\"obius function} $\mu$ and {\it Dirichlet product} $\ast$) which is supported on {\it smooth numbers}. We call $n\in \N$ a $Q-$smooth number by definition when the prime $p|n$ is possible only for $p\le Q$ (to avoid trivialities let $Q\ge 2$). For this, our previous work [C3] inspired us, when looking for multiplicative Ramanujan coefficients (see following definitions) of {\it null-function} $\0(n)\defineq 0$, $\forall n\in\N$. In fact, when adding the hypothesis \lq \lq multiplicative\rq \rq, for these coefficients, say $G:\N \rightarrow \C$, of our $F$, to previous smooth-numbers constraint, we found the results in this paper. Compare section 4 formul\ae, for this kind of heuristic approach. 
\par
Present results hold even for $F\neq \0$, but for null-function we found, so to speak, the best applications. However, for the reason of both giving general results and not going deeply in case $\0$, we will afford the \lq \lq Classification\rq \rq \thinspace for multiplicative Ramanujan coefficients of $\0$ in another, forthcoming paper.  
\medskip
\par
We will give our notation and definitions during the paper, where needed. 
\medskip
\par
\noindent
In this first section, now, we give our Main Theorem, then our Proposition as its \lq \lq local version\rq \rq, and our two Corollaries, as consequences of both previous results; next section contains the statements and proofs of our Lemmas; section 3 gives a brief account of basic hypotheses to get infinite Euler products for general Ramanujan expansions with multiplicative coefficients; section 4 exposes, very quickly, the multiplicative Ramanujan coefficients supported in the smooth numbers; finally, some remarks about forthcoming Classification in the multiplicative part of the \lq \lq cloud\rq \rq (i.e.the Ramanujan coefficients set) of the function $\0$. 
\par
We start giving first definition, as a kind of abbreviation very useful in the sequel.
\par				
\noindent
We say that the arithmetic function $G:\N \rightarrow \C$ {\it is a Ramanujan coefficient}, by definition, when the series 
$$
\sum_{q=1}^{\infty}G(q)c_q(a)
$$
\par
\noindent
converges pointwise in all $a\in \N$. Here $c_q(a)$ is the {\it Ramanujan sum} of $a\in \N$, of modulus $q\in \N$, defined as [M] (compare [R] for the original definition, suggesting why the letter $c$) 
$$
c_q(a)\defineq \sum_{j\le q,(j,q)=1}e^{2\pi i ja/q} 
$$
\par
\noindent
where, as usual, $(j,q)\defineq $g.c.d.$(j,q)$ stands for the {\it greatest common divisor} of these two positive integers. 
\par
In passing, all Ramanujan coefficients $G$ have convergent series (case $a=1$, as \thinspace $c_q(1)=\mu(q)$, $\forall q\in \N$)
$$
\sum_{q=1}^{\infty}G(q)\mu(q). 
$$

\bigskip

\par
\noindent
For any non-empty subset of primes $\P$, say $\corsivoP \subseteq \P$, we define 
$$
\corsivoS_{\corsivoP}\defineq \{ s\in \N : s=1 \; \hbox{\rm or } \; p|s \, \Rightarrow \, p\in \corsivoP\}
$$
\par
\noindent
and for any non-empty subset of natural numbers, say $\corsivoN \subseteq \N$, we define, $\forall r\in \N$, 
$$
(r,\corsivoN)=1\definiz (r,n)=1, \forall n\in \corsivoN. 
$$
\par
\noindent
For $p\in \P$, recall the $p-${\stampatello adic valuation} of $a\in \N$ is \hfill $v_p(a)\defineq \max \thinspace \{ K\in \N_0 : p^K|a \}$, \thinspace  with \thinspace \thinspace $\N_0\defineq \N \cup \{0\}$. 

\bigskip

\par
We start with our central result, about finite Euler products of Ramanujan expansions.
\smallskip 
\par
\noindent {\bf Main Theorem} ({\stampatello finite Euler products}). {\it Let } $G:\N \rightarrow \C$ {\it be multiplicative and choose an arbitrary finite and non-empty subset of primes, say } $\corsivoF \subset \P$, $\corsivoF \neq \emptyset$ ({\it finite}). {\it Assume } $\sum_{(r,\corsivoF)=1}G(r)c_r(a)$ {\it converges pointwise in $a\in \N$. Then, $G$ is a Ramanujan coefficient and the Ramanujan expansion factors as} 
$$
\sum_{q=1}^{\infty}G(q)c_q(a)=\left(\sum_{s\in \corsivoS_{\corsivoF}}G(s)c_s(a)\right) \left(\sum_{(r,\corsivoF)=1}G(r)c_r(a)\right),
\qquad
\forall a\in \N,
$$
\par
\noindent
{\it where the \lq \lq finite factor\rq \rq, at the right hand side, is a }  {\stampatello finite Euler product} : 
$$
\sum_{q=1}^{\infty}G(q)c_q(a)=\left(\prod_{p\in \corsivoF}\sum_{K=0}^{v_p(a)}p^K(G(p^K)-G(p^{K+1}))\right) \left(\sum_{(r,\corsivoF)=1}G(r)c_r(a)\right),
\enspace 
\forall a\in \N.
$$
\par
\noindent
{\it In particular, taking } $G\neq \0$ {\it multiplicative and, for each fixed } $d\in \N$, $\corsivoF$ {\it to be the set of prime divisors of $d$, choosing } $a=1$\enspace {\it this gives, provided the series in RHS converges pointwise in all } $d\in \N$, 
$$
\sum_{q=1}^{\infty}G(q)\mu(q)=\prod_{p|d}(1-G(p))\cdot \left(\sum_{(r,d)=1}G(r)\mu(r)\right), 
\forall d\in \N. 
$$
\smallskip
\par				
\noindent {\bf Proof.} The following formula for the {\bf finite factor}, say, 
$$
\sum_{s\in \corsivoS_{\corsivoF}}G(s)c_s(a)=\prod_{p\in \corsivoF}\sum_{K=0}^{\infty}G(p^K)c_{p^K}(a)=\prod_{p\in \corsivoF}\sum_{K=0}^{v_p(a)}p^K(G(p^K)-G(p^{K+1}))
$$
\par
\noindent
comes from the definition of $\corsivoS_{\corsivoF}$ and the multiplicativity of $G$, then from our Main Lemma (see next $\S2$); so, the finite factor is a finite sum, as a finite product of finite sums (the $p-$Euler factors in Main Lemma). Hence, we may exchange the sum over $s\in \corsivoS_{\corsivoF}$  with the pointwise converging series over $(r,\corsivoF)=1$ as follows: 
$$
\sum_{q=1}^{\infty}G(q)c_q(a)=\sum_{s\in \corsivoS_{\corsivoF}}\sum_{(r,\corsivoF)=1}G(s)c_s(a)G(r)c_r(a)
=\left(\sum_{s\in \corsivoS_{\corsivoF}}G(s)c_s(a)\right) \left(\sum_{(r,\corsivoF)=1}G(r)c_r(a)\right), 
\enspace \forall a\in \N 
$$
\par
\noindent
whence we get $G$ is a Ramanujan coefficient and the first formula. The second one, then, comes from Main Lemma formul\ae. Finally, choosing $\corsivoF$ the prime divisors of $d\in \N$, since now $p\in \corsivoF$ is $p|d$ and $(r,\corsivoF)=1$ is $(r,d)=1$, previous formula gives \enspace $\forall d\in \N$, choosing $a=1$, 
$$
\sum_{q=1}^{\infty}G(q)c_q(1)=\left(\prod_{p|d}\sum_{K=0}^{v_p(1)}p^K(G(p^K)-G(p^{K+1}))\right) \left(\sum_{(r,d)=1}G(r)c_r(1)\right),
$$
\par
\noindent
in which: $v_p(1)=0$ $\forall p\in \P$ and $c_q(1)=\mu(q)$, $c_r(1)=\mu(r)$, whence proving last formula, since $G(1)=1$ (recall: from $G\neq \0$ is multiplicative), and completes the Proof.\hfill $\square$

\medskip

\par
In the following, the series over \thinspace $(r,\corsivoF)=1$ \thinspace will be called the {\bf co-finite factor}.

\medskip

\par
\noindent
We abbreviate 
$$
\P(\corsivoS)\defineq \{ p\in \P : \exists s\in \corsivoS \enspace \hbox{\rm such that} \enspace p|s \}
$$
\par
\noindent
the set of prime factors of all numbers in $\corsivoS \subseteq \N$; trivially, $\P(\{1\})=\emptyset$, while we set by definition $\P(\emptyset)\defineq \emptyset$. In particular, by abuse of notation, we write for the set of prime divisors of $n\in \N$ 
$$
\P(n)\defineq \{ p\in \P : p|n \}. 
$$
\par
Next result is a kind of special case of our Main Theorem, but now \corsivoF \enspace depends on the variable $a\in \N$. 
\smallskip
\par
\noindent {\bf Proposition} ({\stampatello local Euler products}). {\it Let } $G:\N \rightarrow \C$ {\it be multiplicative and assume that the series } $\sum_{(r,a)=1}G(r)\mu(r)$ {\it converges pointwise } $\forall a\in \N$. {\it Then, $G$ is a Ramanujan coefficient and the Ramanujan expansion factors as} 
$$
\forall a\in \N, 
\quad
\sum_{q=1}^{\infty}G(q)c_q(a)=\left(\prod_{p|a}\sum_{K=0}^{v_p(a)}p^K\left(G(p^K)-G(p^{K+1})\right)\right) \left(\sum_{(r,a)=1}G(r)\mu(r)\right). 
$$
\par
\noindent {\bf Proof.} After fixing $a\in \N$, we choose as finite set of primes just $\corsivoF=\P(a)$, so now $p\in \corsivoF$ amounts to $p|a$ and $(r,\corsivoF)=1$ means $(r,a)=1$, whence $c_r(a)=\mu(r)$ : everything else is like in Main Theorem Proof.\hfill $\square$ 

\bigskip

\par
\noindent
In order to give the main consequences of our Main Theorem, we have to define the following sets: 
$$
\corsivoF(G)\defineq \{ p\in \P : G(p)=1 \}
$$
\par
\noindent
and its subset 
$$
\corsivoFzero(G)\defineq \{ p\in \P : G(p^K)=1, \forall K\in \N_0 \}
$$
\par
\noindent
for which sets we explicitly point out the finiteness, when $G$ is a Ramanujan coefficient, as follows. 
\smallskip
\par				
\noindent {\bf Remark 1.} In case $G$ is a Ramanujan coefficient, in particular, at $a=1$ the series
$$
\sum_{q=1}^{\infty}G(q)\mu(q) \quad \hbox{\rm converges}.
$$
\par
\noindent
This entails, from the necessary condition for series to converge, that 
$$
G(q)\to 0,
$$
\par
\noindent
as $q\to \infty$ in the set of square-free numbers $q$; in particular, 
$$
G(p)\to 0,
$$
\par
\noindent
as $p\to \infty$ in the set of prime numbers $\P$. This, in turn, implies that, for any Ramanujan coefficient $G:\N \rightarrow \C$, the set $\corsivoF(G)$, whence $\corsivoFzero(G)\subseteq \corsivoF(G)$ too, are finite sets.\hfill $\diamond$

\bigskip

\par
We give our first application to the cloud ($=$set of Ramanujan coefficients [C2],[C4]) of null-function $\0$. 
\smallskip
\par
\noindent {\bf Corollary 1.} {\it Let } $G:\N \rightarrow \C$ {\it be } {\stampatello multiplicative}, {\it with } $
\corsivoF(G) = \emptyset$ {\it and assume that } $\sum_{(r,a)=1}G(r)\mu(r)$ {\it converges pointwise } $\forall a\in \N$. {\it Then, $G$ is a }  {\stampatello Ramanujan coefficient} {\it and} 
$$
\sum_{q=1}^{\infty}G(q)c_q(a)=\0(a)
\enspace \Longleftrightarrow \enspace 
\sum_{q=1}^{\infty}G(q)\mu(q)=0. 
$$
\par
\noindent {\bf Proof.} The implication $\Rightarrow$ follows trivially from $a=1$. Notice the convergence in RHS. 
\par
\noindent
In order to prove $\Leftarrow$ we apply our Proposition: 
$$
\sum_{q=1}^{\infty}G(q)c_q(a)=\left(\prod_{p|a}\sum_{K=0}^{v_p(a)}p^K (G(p^K)-G(p^{K+1}))\right)\left(\sum_{(r,a)=1}G(r)\mu(r)\right)
=0, 
\quad 
\forall a\in \N, 
$$
\par
\noindent
since our Lemma called \lq \lq Fact 2\rq \rq, in next section $\S2$, gives in particular 
$$
\sum_{(r,a)=1}G(r)\mu(r)
\enspace {\rm converges} \enspace {\rm pointwise} \enspace  
\forall a\in \N 
\enspace {\rm and} \enspace 
\sum_{q=1}^{\infty}G(q)\mu(q)=0
\quad \Longrightarrow \quad 
\sum_{(r,a)=1}G(r)\mu(r)=0, 
\enspace
\forall a\in \N. 
$$
\par
\noindent
Thus, our task is completed.\hfill $\square$

\bigskip

\par
We give our second application, even easier to prove, to the cloud of $\0$. 
\smallskip
\par
\noindent {\bf Corollary 2.} {\it Let } $G:\N \rightarrow \C$ {\it be } {\stampatello multiplicative}, {\it with } $
\corsivoFzero(G) \neq \emptyset$ {\it and assume that } $\sum_{(r,\corsivoFzero(G))=1}G(r)c_r(a)$ {\it converges pointwise } $\forall a\in \N$. {\it Then, $G$ is a } {\stampatello Ramanujan coefficient} {\it and} 
$$
\sum_{q=1}^{\infty}G(q)c_q(a)=\left(\prod_{p\in \corsivoFzero(G)}\sum_{K=0}^{v_p(a)}p^K(G(p^K)-G(p^K+1))\right)\left(\sum_{(r,\corsivoFzero(G))=1}G(r)c_r(a)\right)
=\0(a). 
$$
\par
\noindent {\bf Proof.} Since \enspace $\corsivoFzero(G) \neq \emptyset$, pick up a prime \enspace $p \in \corsivoFzero(G)$ \enspace and consider its Euler factor: 
$$
\sum_{K=0}^{\infty}G(p^K)c_{p^K}(a)=\sum_{K=0}^{v_p(a)}p^K\left(G(p^K)-G(p^{K+1})\right)=\0(a), 
$$
\par
\noindent
from our Main Lemma in $\S2$; our Main Theorem completes the proof, choosing \enspace $\corsivoF = \corsivoFzero(G)$.\hfill $\square$ 

\medskip

\par
\noindent {\bf Remark 2.} We call {\stampatello standard} the factorization above, for $G$ in the hypotheses above. 
\par
\noindent
Hereafter, when not specified, equations depending on $a$ hold \enspace $\forall a\in \N$. 
\par
\noindent
We explicitly point out that this co-finite factor may vanish for {\stampatello some} $a\in \N$, but {\stampatello not all} of them.\hfill $\diamond$

\vfill
\eject

\par				
\noindent {\bf 2. Lemmata}

\bigskip

\par
\noindent
In the following, the symbol QED (Quod Erat Demonstrandum$=$what was to be shown) will delimit small parts of a Proof; in following sections, it will  even be the end of a small result's proof. 

\medskip

\par
We start with a very short Lemma, a \lq \lq {\stampatello fact}\rq \rq, allowing a quicker proof of next Main Lemma. 
\smallskip
\par
\noindent {\bf Fact 1.} {\it For all } $p\in \P$, $K\in \N_0$ {\it and } $a\in \N$ {\it we have}: 
$$
c_{p^K}(a)=c_{p^K}(p^{v_p(a)})=\varphi(p^K)\cdot {{\mu(p^{K-\min(K,v_p(a))})}\over {\varphi(p^{K-\min(K,v_p(a))})}}. 
$$
\par
\noindent {\bf Proof.} Assuming $K\in \N$ henceforth, as case $K=0$ is trivially true (recall $c_1(a)=1$, for all $a\in \N$), first equation follows from writing $a=mp^{v_p(a)}$, with $(m,p)=1$, and the definition of Ramanujan sum [R] (compare [M], too) that we recall, for the imaginary exponential $e_q(n)\defineq e^{2\pi in/q}$, 
$$
c_{p^K}(a)\defineq \sum_{j\le p^K, (j,p)=1}e_{p^K}(ja)
=\sum_{j\le p^K, (j,p)=1}e_{p^K}(jmp^{v_p(a)})
=\sum_{j'\le p^K, (j',p)=1}e_{p^K}(j'p^{v_p(a)})
=c_{p^K}(p^{v_p(a)}), 
$$
\par
\noindent
from the invertible change of variables $j'\equiv jm(\bmod p^K)$; while second equation follows from the 1936 H\"older formula, see page 22 [CM] (and [D], page 149, for a proof) : 
$$
c_q(a)=\varphi(q)\cdot {{\mu(q/(q,a))}\over {\varphi(q/(q,a))}}
\enspace 
\forall q\in \N, \forall a\in \N, 
$$
\par
\noindent
since, of course, the greatest common divisor is \enspace $(p^K,p^{v_p(a)})=p^{\min(K,v_p(a))}$.\hfill $\square$

\medskip

\par
The core of our exposition is the following main result. 
\smallskip
\par
\noindent {\bf Main Lemma} ({\stampatello $p-$Euler factors}). {\it Let } $G:\N \rightarrow \C$ {\it be any arithmetic function and fix any } $p\in \P$. {\it Then, } $\forall a\in \N$, 
$$
\sum_{K=0}^{\infty}G(p^K)c_{p^K}(a)=\sum_{K=0}^{v_p(a)}G(p^K)\varphi(p^K)-G(p^{v_p(a)+1})p^{v_p(a)}, 
$$
\par
\noindent
{\it from which}
$$
\sum_{K=0}^{\infty}G(p^K)c_{p^K}(a)=\sum_{K=0}^{v_p(a)}p^K\left(G(p^K)-G(p^{K+1})\right), 
$$
\par
\noindent
{\it just reordering $p-$th powers.}
\par
\noindent
{\it Furthermore,}
$$
\sum_{K=0}^{\infty}G(p^K)c_{p^K}(a)=\0(a)
\quad \Longleftrightarrow \quad
G(p^K)=G(1), \enspace \forall K\in \N. 
$$
\smallskip
\par
\noindent {\bf Proof.} From Fact 1, we easily obtain 
$$
c_{p^K}(a)=\varphi(p^K), 
\enspace {\rm if} \enspace
0\le K\le v_p(a), 
$$
$$
c_{p^K}(a)=\varphi(p^{v_p(a)+1})\cdot {{\mu(p)}\over {\varphi(p)}}
=-(p^{v_p(a)+1}-p^{v_p(a)})/(p-1)
=-p^{v_p(a)}, 
\enspace {\rm if} \enspace
K=v_p(a)+1 
$$
\par
\noindent
and 
$$
c_{p^K}(a)=0, 
\enspace \forall K>v_p(a)+1.  
$$
\par
\noindent
These formul\ae, of course, give instantly our first equation. \hfill QED(1st eq.)
\par
The second equation comes from first equation, since 
$$
\varphi(p^K)=p^K-p^{K-1},
\enspace 
\forall K\in \N
$$
\par				
\noindent
and, trivially, $\varphi(p^0)=1$, then we rearrange $p-$th powers, with Abel summation [T] trick (whose details are left to the reader). \hfill QED(2nd eq.)
\par
For the equivalence, notice that implication \thinspace \lq \lq $\Leftarrow$\rq \rq \thinspace immediately follows from our second equation, proved just now. \hfill QED($\Leftarrow$)
\par
Proving \thinspace \lq \lq $\Rightarrow$\rq \rq \thinspace is all we need to show, to conclude. 
\par
\noindent
Our starting hypothesis, now, is that the Euler $p-$factor with coefficients $G$ is $0$ in all $a\in \N$, which is equivalent to the following, say: 
$$
\sum_{K=0}^{\infty}G(p^K)c_{p^K}(a)=0, \enspace \forall a\in \N
\quad \Longleftrightarrow \quad
E_{p,G}(v)\defineq \sum_{K=0}^{v}p^K (G(p^K)-G(p^{K+1}))=0, \enspace \forall v\in \N_0,
\leqno{(\ast)}
$$
\par
\noindent
again from our second equation above. We may solve, in the definition in RHS of $(\ast)$, for $G(p^{v+1})$, getting 
$$
G(p^{v+1})=G(p^v)-p^{-v} (E_{p,G}(v)-E_{p,G}(v-1)), \enspace \forall v\in \N
$$
\par
\noindent
and this implies, under hypothesis $(\ast)$,
$$
G(p^{v+1})=G(p^v), \enspace \forall v\in \N
\quad \Longrightarrow \quad 
G(p^K)=G(1), \enspace \forall K\in \N, 
$$
\par
\noindent
by induction on $v\in \N$ (case $v=0$ is trivial). \hfill QED($\Rightarrow$) 
\par
\noindent
The proof is complete.\hfill $\square$ 

\bigskip

\par
We will give the other Lemma, again a \lq \lq {\stampatello fact}\rq \rq, as it is a straight consequence of our Main Theorem; which, in turn, is proved applying previous Lemma. 
\smallskip
\par
\noindent {\bf Fact 2.} {\it If } $G:\N \rightarrow \C$ {\it is multiplicative}, $\corsivoF(G)=\emptyset$ {\it and } $\sum_{(r,d)=1}G(r)\mu(r)$ {\it converges pointwise } $\forall d\in \N$, {\it then} 
$$
\sum_{(r,d)=1}G(r)\mu(r)=\0(d)
\enspace \Longleftrightarrow \enspace
\sum_{q=1}^{\infty}G(q)\mu(q)=0. 
$$
\smallskip
\par
\noindent {\bf Proof.} Since the \lq \lq $\Longrightarrow$\rq \rq \enspace part follows trivially for $d=1$, we have only to prove the \lq \lq $\Longleftarrow$\rq \rq \enspace : from last part of our Main Theorem (see $\S1$) we get 
$$
\sum_{q=1}^{\infty}G(q)\mu(q)=\prod_{p|d}(1-G(p))\cdot \sum_{(r,d)=1}G(r)\mu(r),
\quad
\forall d\in \N
$$
\par
\noindent
and the hypothesis on \enspace $\corsivoF(G)$ \enspace implies 
$$
\prod_{p|d}(1-G(p))\neq 0,
\quad
\forall d\in \N, 
$$
\par
\noindent
whence in two lines we get the \lq \lq $\Longleftarrow$\rq \rq.\hfill $\square$ 

\vfill
\eject

\par				
\noindent {\bf 3. A glimpse into infinite Euler products for Ramanujan expansions}

\bigskip

\par
\noindent
We start, for infinite Euler products, from a {\bf new definition}. 
\par
We say a Ramanujan expansion {\it is an Euler product} iff (if and only if) $G:\N \rightarrow \C$ is a multiplicative Ramanujan coefficient and, for all fixed $a\in \N$, we have 
$$
\lim_Q \sum_{q\le Q}G(q)c_q(a) = \lim_Q \sum_{{q=1}\atop {p|q\Rightarrow p\le Q}}^{\infty}G(q)c_q(a)
\leqno{(1)}
$$
\par
\noindent
(notice: LHS is the Ramanujan expansion itself), whence in one line we get the infinite Euler product: 
$$
\sum_{q=1}^{\infty}G(q)c_q(a) = \lim_Q \prod_{p\le Q}\sum_{K=0}^{\infty}G(p^K)c_{p^K}(a) 
= \prod_{p\in \P} \sum_{K=0}^{\infty}G(p^K)c_{p^K}(a), 
$$
\par
\noindent
since the second equation is itself the {\it definition} of {\bf infinite Euler product} for a Ramanujan expansion. 
\par
\noindent
(Here we do not say if the product vanishes or not: compare Property 2 \& Remark 3, following.)
\medskip
\par
Notice in $(1)$ the link with $Q-$smooth numbers, in its RHS : see $\S4$. Compare also [C1]. 
\medskip
\par
\noindent
The {\bf absolute convergence of} the {\bf Ramanujan expansion} (since $G(q)\in \C$, following $|\enspace |$ are moduli) : 
$$
\forall a\in \N,
\qquad
\sum_{q=1}^{\infty}\left| G(q)c_q(a)\right| \enspace < \thinspace \infty 
\leqno{(2)}
$$
\par
\noindent
implies : it's an Euler product. We prove it, even if it's a classic in the literature, as following little result. 
\smallskip
\par
\noindent {\bf Property 1.} {\it If a Ramanujan expansion satisfies $(2)$, with multiplicative } $G$, {\it then it is an Euler product}. 
\smallskip
\par
\noindent {\bf Proof.} Given \thinspace $a\in \N$ \thinspace and \thinspace $\varepsilon>0$, there exists a \enspace $Q=Q(a,\varepsilon)\in \N$ \enspace such that (from $(2)$ hypothesis) 
$$
\sum_{q>Q}\left| G(q)c_q(a)\right|<\varepsilon, 
$$
\par
\noindent
implying: (since the $q\le Q$ have, trivially, all prime divisors $p\le Q$)
$$
\left| \prod_{p\le Q}\sum_{K=0}^{\infty}G(p^K)c_{p^K}(a) - \sum_{q\le Q}G(q)c_q(a)\right| \le 
 \sum_{{q>Q}\atop {p|q \Rightarrow p\le Q}}\left| G(q)c_q(a)\right|
\le \sum_{q>Q}\left| G(q)c_q(a)\right|
< \varepsilon, 
$$
\par
\noindent
i.e., the limit over $Q\to \infty$ proves $(1)$, so the Ramanujan expansion is an Euler product.\hfill QED 
\medskip
\par
Now, of course, we may generalize the results we obtained above for the product of a finite factor and a co-finite factor (for Ramanujan expansions) to \lq \lq arbitrary products\rq \rq, of two factors; but over two subsets, say $\corsivoS$ \enspace and $\corsivoR$, that are \lq \lq coprime\rq \rq, whose \lq \lq product\rq \rq, say a {\stampatello coprime product}, is the set of natural numbers. We do it for only two \lq \lq factors\rq \rq, but this generalizes to a finite number of factors (for $\N$ here).
\par
We say that \thinspace $\corsivoS$ \thinspace and \thinspace $\corsivoR$ \thinspace are {\bf coprime factors} of natural numbers, writing $\corsivoS  \otimes \corsivoR = \N$, where $\otimes$ is the {\bf coprime product}, iff \enspace $\forall a\in \N$, $\exists ! (r,s)\in \corsivoR \times \corsivoS$, (a unique couple in the cartesian product), with $r$ and $s$ coprime each other, such that : $a=rs$. These coprime factors in $\corsivoR \otimes \corsivoS = \N$ have a particular shape, in the family of $\N$ subsets, characterized by the set of prime numbers dividing all the naturals in the set. In fact, 
$$
\corsivoS \cong \corsivoP
$$
\par
\noindent
with $\cong$ a bijective function from the family of $\N$ subsets, of the form $\corsivoS = \{ s\in \N : s=1 \enspace \hbox{\rm or} \enspace p|s \Rightarrow p\in \corsivoQ \}$ (for a fixed $\corsivoQ \subseteq \P$), to the family of subsets of primes $\P$, given by both
$$
\corsivoP = \P(\corsivoS)
\quad \hbox{\rm and\enspace its\enspace inverse} \quad
\corsivoS = \corsivoS_{\corsivoP}
$$
\par				
\noindent
(see the definitions in $\S1$), which also give the link, writing $\oplus$ as usual for the {\it disjoint union} of sets,  
$$
\corsivoS \otimes \corsivoR = \N
\quad \Longleftrightarrow \quad
\P(\corsivoS) \oplus \P(\corsivoR) = \P
$$
\par
\noindent
The trivial case \enspace $\{1\}\otimes \N=\N\otimes \{1\}=\N$ \enspace corresponds to: $\P(1)\oplus \P(\N)=\P(\N)\oplus \P(1)=\P$. 

\medskip

\par
We avoid the trivial case (so we'll write \lq \lq non-trivial\rq \rq, in the following). 

\medskip

\par
\noindent
The (non-trivial) coprime factorization $\N=\corsivoS \otimes \corsivoR$ gives for absolutely converging Ramanujan expansions (see $(2)$ above), $\forall a\in \N$, convergence of both following series and 
$$
\sum_{q=1}^{\infty}G(q)c_q(a)=\left(\sum_{r\in \corsivoR}G(r)c_r(a)\right)\left(\sum_{s\in \corsivoS}G(s)c_s(a)\right) 
$$
\par
\noindent
whence, above Property 1 implies the infinite Euler product
$$
\sum_{q=1}^{\infty}G(q)c_q(a)=\prod_{p\in \P(\corsivoR)}\sum_{K=0}^{\infty}G(p^K)c_{p^K}(a) \cdot \prod_{p\in \P(\corsivoS)}\sum_{K=0}^{\infty}G(p^K)c_{p^K}(a), 
\leqno{(3)}
$$
\par
\noindent
which is, of course, a generalization of our Main Theorem for any couple $\corsivoR \otimes \corsivoS = \N$, without the constraint of having one of the two factors finite (namely, $|\P(\corsivoS)|<\infty$, see $\S1$). 
\par
This is possible (see Property 1 above), of course, assuming the absolute convergence for the Ramanujan expansion (with multiplicative coefficients).
\par
\noindent
However a \lq \lq side effect\rq \rq, so to speak, of this assumption is that the infinite Euler products may not vanish. Also this is a classic and known result that we will prove very briefly, in our next Property 2, following. 
\par
The \lq \lq local Euler product\rq \rq, compare Proposition $\S1$, follows from $\corsivoS = \P(a)$ and $\corsivoR = \{ r\in \N : (r,a)=1\}$ in $(3)$, so \enspace $\forall a\in \N$, 
$$
\sum_{q=1}^{\infty}G(q)c_q(a)=\prod_{p\not \thinspace \vert \thinspace a}(1-G(p)) \cdot \prod_{p|a}\sum_{K=0}^{\infty}G(p^K)c_{p^K}(a). 
$$
\par
\noindent
Notice : the (infinite) Euler product for the co-finite factor in our Proposition (see $\S1$) needs the absolute convergence $(2)$, in order to get from $(1)$ following equation 
$$
\forall a\in \N, 
\quad 
\sum_{(r,a)=1}G(r)\mu(r)=\prod_{p\not \thinspace \vert \thinspace a}(1-G(p)). 
$$
\par
\noindent
The absolute convergence (see $(2)$ above) implies, as we prove now (compare $\S1$ beginning,Chapter II [K]), 
$$
\prod_{p\not \thinspace \vert \thinspace a}(1-G(p))=\prod_{p\notin \P(a)}(1-G(p))
\neq 0. 
$$
\par
Actually, since $\P(a)$ is a finite subset of primes, we'll prove this for any finite (non-empty) $\corsivoF \subset \P$. 
\smallskip
\par
\noindent {\bf Property 2.} {\it Let } $G:\N \rightarrow \C$ {\it be a multiplicative Ramanujan coefficient of an absolutely convergent Ramanujan expansion. Then, for any finite subset of primes } $\corsivoF \neq \emptyset$, {\it the Euler product} 
$$
\prod_{p\notin \corsivoF}(1-G(p))
$$
\par
\noindent
{\it converges to a non-zero complex number, provided } $G(p)=1$ $\Rightarrow$ $p\in \corsivoF$.
\par
\noindent {\bf Proof.} Using \enspace $p\notin \corsivoF$ $\Rightarrow$ $G(p)\neq 1$, as a consequence of Theorem 6 in $\S2.2$, Chapter 5 of [A], we get 
$$
\sum_{p\notin \corsivoF}|G(p)|<\infty 
\quad \Longrightarrow \quad 
\prod_{p\notin \corsivoF}(1-G(p))\neq 0 
$$
\par				
\noindent
(implicit: the infinite product converges) and the Ramanujan expansion absolute convergence at $a=1$ (compare Remark 1): 
$$
\sum_{q=1}^{\infty}\mu^2(q)|G(q)|<\infty 
$$
\par
\noindent
implies: 
$$
\sum_{p\in \P}|G(p)|<\infty, 
$$
\par
\noindent
whence 
$$
\sum_{p\notin \corsivoF}|G(p)|<\infty 
$$
\par
\noindent
concludes.\hfill QED 
\medskip
\par
\noindent {\bf Remark 3.} The same hypotheses of Property 2 give, provided \enspace $\corsivoF(G)=\emptyset$, 
$$
\prod_p (1-G(p))\neq 0,
$$
\par
\noindent
as it is, from previous proof, crystal clear.\hfill $\diamond$ 
\medskip
\par
Then, this Remark implies that Corollary 1 with the additional hypothesis $(2)$ (of absolute convergence) would supply a $G$ NOT in the cloud of $\0$; in fact, adding $(2)$, Property 1 gives the infinite product: 
$$
\sum_{q=1}^{\infty}G(q)\mu(q)=\prod_p (1-G(p))=0, 
$$
\par
\noindent
impossible from Remark 3, so $G$ is not in $\0-$cloud. Compare Corollary 3 in $\S5$, that's linked to Property 2. 
\medskip
\par
However, since Corollary 2 holds without conditions, on these kind of infinite Euler products, we may have the additional hypothesis $(2)$ in this Corollary and still get, surprisingly enough, a $G:\N \rightarrow \C$ in the cloud of $\0$. In other words, there exist multiplicative Ramanujan coefficients $G:\N \rightarrow \C$ with $\corsivoFzero(G)\neq \emptyset$ and absolutely convergent Ramanujan expansion, namely satisfying $(2)$ above, that are in $\0$ cloud. As an example, take $G(q)\defineq 1/q^3$ on odd numbers and $G(2^K)\defineq 1$, $\forall K\in \N_0$ : it has $\corsivoFzero(G)=\{2\}$, so $G$ is in $\0-$cloud (Corollary 2) and the trivial bound $|c_q(a)|\le \varphi(q)\le q$, on all $a\in \N$, gives $(2)$. In fact: (next $K-$series is a finite sum depending on $a$, after Main Lemma, and compare Remark 2) 
$$
\sum_{q=1}^{\infty}|G(q)c_q(a)|=\sum_{K=0}^{\infty}|c_{2^K}(a)|\cdot \sum_{{q=1}\atop {q\not \equiv 0(\!\!\bmod 2)}}^{\infty}{{|c_q(a)|}\over {q^3}}
\ll_a \sum_{{q=1}\atop {q\not \equiv 0(\!\!\bmod 2)}}^{\infty}{{\varphi(q)}\over {q^3}}
\ll_a \sum_{q=1}^{\infty}{1\over {q^2}}
\ll_a 1,
\enspace 
\forall a\in \N. 
$$

\bigskip

\par
\noindent {\bf 4. A link to some Ramanujan coefficients supported on smooth numbers}

\bigskip

\par
\noindent
See that, if $G:\N \rightarrow \C$ is a {\stampatello multiplicative Ramanujan coefficient}, then the series in RHS of $(1)$ of $\S3$ becomes, \enspace $\forall a \in \N$ :  
$$
\sum_{{q=1}\atop {p|q\Rightarrow p\le Q}}^{\infty}G(q)c_q(a)=\prod_{p\le Q}\sum_{K=0}^{\infty}G(p^K)c_{p^K}(a)
=\prod_{p\le Q}\sum_{K=0}^{v_p(a)}p^K (G(p^K)-G(p^{K+1})), 
$$
\par
\noindent
again, by our Main Lemma, so the LHS is a {\bf finite Ramanujan expansion} (see [CMS], compare [CM] $\S4$) but, in the terminology of [CM], it is {\bf not} a {\bf pure} one, as it depends on $a\in \N$, since RHS above has inside $v_p(a)$ this dependence. 
\par
In other words, we have just proved the following small result. 
\smallskip
\par				
\noindent {\bf Property 3.} {\it Let } $G:\N \rightarrow \C$ {\it be a multiplicative Ramanujan coefficient, supported on the } $Q-${\it smooth numbers } ({\it namely}, $\exists p>Q$ : $p|q$ $\Rightarrow $ $G(q)=0$). {\it Then the Ramanujan expansion is } {\stampatello finite}. 
\medskip
\par
Here $G$ support, abbrev. $\supporto(G)\defineq \{ q\in \N : G(q)\neq 0\}$, is included in the $Q-$smooth numbers, which make up an infinite set (for each fixed $Q\ge 2$, of course). However, the expansion becomes finite as an immediate effect of second formula in our Main Lemma (which also explains why, see the above, the Ramanujan expansion {\it length} depends on $a$). Our papers [C2],[C4] characterize the finite \& pure Ramanujan expansions : they are truncated divisor sums. This is not the case, here, because \lq \lq finite and pure\rq \rq, in particular, implies that the expansion has a fixed length, not dependent on the variable $a$. 
\medskip
\par
Compare, also, the Ramanujan expansion we gave in [C3] for all $F$ satisfying the Ramanujan Conjecture and with Eratosthenes transform $F'$ supported in $Q-$smooth numbers: see $\S1$ beginning. A result, in [C3], also allows to give an explicit formula for Ramanujan coefficients, named after Wintner, i.e. $G=\Wintner(F)$ (see Theorem 1 there). These have the property : $\supporto(F')$ is inside $Q-$smooth n.s $\Rightarrow$ $\supporto(\Wintner(F))$ is inside $Q-$smooth n.s; so, if we know that $F'$ vanishes outside $Q-$smooth numbers, then : its Ramanujan expansion with Wintner coefficients not only converges to $F$ (by quoted Theorem 1), but, also, these Ramanujan coefficients are supported on $Q-$smooth numbers. However, in order to apply above Property 3 we should also know if: $\Wintner(F)$ is a multiplicative function or not, that, in general, is not known! 

\bigskip

\par
\noindent {\bf 5. A coming soon for a classification of multiplicative Ramanujan coefficients in $\0-$cloud}

\bigskip

\par
\noindent
The third possibility for multiplicative Ramanujan coefficients, apart from that of Corollary 1 (i.e., $\corsivoF(G)=\emptyset$) and that of Corollary 2 (i.e., $\corsivoFzero(G)\neq \emptyset$) is the case: $\corsivoF(G)\neq \emptyset$ and $\corsivoFzero(G)=\emptyset$, in our third Corollary. 
\smallskip
\par
\noindent {\bf Corollary 3.} {\it Let } $G:\N \rightarrow \C$ {\it be } {\stampatello multiplicative}, {\it with } $\corsivoF(G)\neq \emptyset$, $\corsivoFzero(G)=\emptyset$ {\it and assume that } $\sum_{(r,a)=1}G(r)\mu(r)$ {\it converges pointwise } $\forall a\in \N$. {\it Then $G$ is a }  {\stampatello Ramanujan coefficient} {\it and} 
$$
\sum_{q=1}^{\infty}G(q)c_q(a)=\0(a) 
\quad 
\Longleftrightarrow 
\quad 
\sum_{(r,\corsivoF(G))=1}G(r)\mu(r)=0. 
$$
\medskip
\par
Thus, just alike Corollary 1 can't have the additional hypothesis of absolute convergence $(2)$ above (see $\S3$, soon after Remark 3), our Property 2 in $\S3$ gives the same for Corollary 3; in fact, under $(2)$ Property 1 gives an infinite product: 
$$
\sum_{(r,\corsivoF(G))=1}G(r)\mu(r)=\prod_{p\notin \corsivoF(G)}(1-G(p))=0
$$
\par
\noindent
which can not vanish by Property 2, so $G$ is not in $\0-$cloud. 
\medskip
\par
As we saw above (in $\S3$, again after Remark 3), in Corollary 2 we may assume absolute convergence $(2)$; now, we know that it's the only case in which we may do it (Corollaries 1 and 3 with $(2)$, so to speak, go out of $\0-$cloud). Since both Ramanujan's and Hardy's coefficients are in the scope of Corollary 1, for Ramanujan's and of Corollary 3, for Hardy's (as $G(q)=1/\varphi(q)$ has $G(p)=1$ iff $p=2$, Hardy's have $\corsivoF(G)=\{2\}$), we have a theoretical, general reason for them not to satisfy, of course, $(2)$ (absolute convergence above). 

\medskip

\par
This Corollary completes the {\bf Classification of multiplicative Ramanujan coefficients of } $\0$. 
\par
\noindent
In a forthcoming paper, we'll prove this last and more difficult case (Corollaries 1 \& 2 prove other two cases). 

\vfill

\par
\noindent {\bf Acknowledgments.} I wish to thank once again Ram Murty for our papers together, a never ending source of inspiration in the world of Ramanujan expansions. Also, I wish to thank Luca Ghidelli for a simplification, in the Proof of Main Lemma. Last but not least, I wish to thank both Maurizio Laporta and Luca Ghidelli for their comments improving the exposition. 
\smallskip
\par
By the way,I invite all the readers to give their feedback as soon as possible. 

\eject

\par				
\centerline{\stampatello Bibliography}

\bigskip

\item{[A]} L.V. Ahlfors, {\sl Complex Analysis}, An introduction to the theory of analytic functions of one complex variable. Third edition. International Series in Pure and Applied Mathematics. McGraw-Hill Book Co., New York, 1978. 
\smallskip
\item{[C1]} G. Coppola, {\sl An elementary property of correlations}, Hardy-Ramanujan J. {\bf 41} (2018), 68--76. Available online 
\smallskip
\item{[C2]} G. Coppola, {\sl A map of Ramanujan expansions}, ArXiV:1712.02970v2. (Second Version) 
\smallskip
\item{[C3]} G. Coppola,{\sl A smooth shift approach for a Ramanujan expansions},ArXiV:1901.01584v3.(Third Version) 
\smallskip
\item{[C4]} G. Coppola, {\sl Recent results on Ramanujan expansions with applications to correlations}, to appear on Rend. Semin. Mat. Univ. Politec. Torino 
\smallskip
\item{[CMS]} G. Coppola, M. Ram Murty and B. Saha, {\sl Finite Ramanujan expansions and shifted convolution sums of arithmetical functions},  J. Number Theory {\bf 174} (2017), 78--92. 
\smallskip
\item{[CM]} G. Coppola and M. Ram Murty, {\sl Finite Ramanujan expansions and shifted convolution sums of arithmetical functions, II}, J. Number Theory {\bf 185} (2018), 16--47. 
\smallskip
\item{[D]} H. Davenport, {\sl Multiplicative Number Theory}, 3rd ed., GTM 74, Springer, New York, 2000. 
\smallskip
\item{[K]} A. Karatsuba, {\sl Basic Analytic Number Theory}, Translated from the second (1983) Russian edition and with a preface by Melvyn B. Nathanson. Springer-Verlag, Berlin, 1993. 
\smallskip
\item{[M]} M. Ram Murty, {\sl Ramanujan series for arithmetical functions}, Hardy-Ramanujan J. {\bf 36} (2013), 21--33. Available online 
\smallskip
\item{[R]} S. Ramanujan, {\sl On certain trigonometrical sums and their application to the theory of numbers}, Transactions Cambr. Phil. Soc. {\bf 22} (1918), 259--276. 
\smallskip
\item{[T]} G. Tenenbaum, {\sl Introduction to Analytic and Probabilistic Number Theory}, Cambridge Studies in Advanced Mathematics, {46}, Cambridge University Press, 1995. 
\smallskip
\item{[W]} A. Wintner, {\sl Eratosthenian averages}, Waverly Press, Baltimore, MD, 1943. 

\vfill

\leftline{\tt Giovanni Coppola - Universit\`{a} degli Studi di Salerno (affiliation)}
\leftline{\tt Home address : Via Partenio 12 - 83100, Avellino (AV) - ITALY}
\leftline{\tt e-mail : giovanni.coppola@unina.it}
\leftline{\tt e-page : www.giovannicoppola.name}
\leftline{\tt e-site : www.researchgate.net}

\bye